\documentclass{IOS-Book-Article}

\usepackage{ amsmath,amssymb} 
\usepackage{amsthm}

\usepackage{hyperref}

\hypersetup{
  colorlinks   = true, 
  urlcolor     = blue, 
  linkcolor    = blue, 
  citecolor   = red 
}

\usepackage[msc-links]{amsrefs}

\newtheorem{thm}{Theorem}

\newtheorem {prop}{Proposition}
\newtheorem {lem}{Lemma}

\newtheorem {exa}{Example}
\newtheorem {exe}{Exercise}

\newtheorem {rem}{Remark}
\newtheorem {cor}{Corollary}

\newcommand\Z{\mathbb Z}
\newcommand\Q{\mathbb Q}
\newcommand\R{\mathbb R}
\newcommand\C{\mathbb C}
\def\P{\mathbb P}
\def\H{\mathcal H}
\def\L{\mathcal L}

\newcommand\M{\mathcal M}
\newcommand\A{\mathcal A}

\newcommand\X{\mathcal X}     

\newcommand\E{\mathcal E}

\newcommand\HS{\mathfrak H}

\newcommand\m{\mathfrak m}
\newcommand\an{\mathfrak a}
\newcommand\bn{\mathfrak b}
\newcommand\en{\mathfrak e}
\newcommand\hn{\mathfrak h}

\newcommand\p{\mathfrak p}
\newcommand\f{{\mathfrak a}}
\newcommand\h{\mathfrak b}
\newcommand\g{\mathfrak c}

\newcommand\<{\langle}
\def\>{\rangle}   

\def\a{\alpha}   
\def\b{\beta}
\def\t{\tau}
\newcommand\s{\sigma}

\def\T{\theta} 
\def\th{\theta}
\def\d{\delta }

\def\O{\Omega}
\newcommand\om{\omega}

\newcommand\G{\Gamma}

\newcommand\xs{{\Bbb o}}

\newcommand\iso{{\, \cong\, }}

\newcommand\sem{\rtimes}

\newcommand\embd{\hookrightarrow}

\newcommand\Aut{\mbox{Aut}}
\newcommand\J{\mbox{Jac }}
\newcommand\Jac{\mbox{Jac }}

\newcommand\e{\eta} 
\def\e{\varepsilon}

\newcommand\bC{\textbf{C}}

\newcommand\JJ{G}
\newcommand\rmJ{H}

\def\sn{\mbox{sn }}
\def\cn{\mbox{cn }}
\def\dn{\mbox{dn }}

\newcommand{\ch}[2]
{\begin{bmatrix}
 #1 \\
 #2\\
\end{bmatrix}}

\begin{document}
\begin{frontmatter}          
%
\title{Theta functions  of superelliptic curves}

\runningtitle{Thetanulls of superelliptic curves}


\author[A]{\fnms{Lubjana} \snm{Beshaj}} 
\author[B]{\fnms{Artur} \snm{Elezi}} 
\author[C]{\fnms{Tony} \snm{Shaska}}

\runningauthor{Beshaj/Elezi/Shaska}

\address[A]{Oakland University, \\Rochester, MI, USA; \\ E-mail: beshaj@oakland.edu}
\address[B]{American University,  \\Washington, DC, USA;   \\E-mail: aelezi@american.edu}
\address[C]{Oakland University, \\Rochester, MI, USA;  \\E-mail: shaska@oakland.edu}

\begin{abstract}
In this short survey we give a description of the theta functions of algebraic curves, half-integer theta-nulls, and the fundamental theta functions.  
We describe how to determine such fundamental theta functions and describe the components of the moduli space in terms of such functions. Several open problems are suggested. 
\end{abstract}


\end{frontmatter}


\maketitle


\section{Introduction}

Superelliptic curves are special in many ways in algebraic geometry and number theory.  The point of this volume is to emphasize such special points as illustrated in \cite{nato-2} and throughout other papers of this volume.  One of the most special properties of superelliptic curves (i.e., which is not known to be true for non-superelliptic curves) are the properties of the theta functions.  In other words, for the experts this means that the Thomae's formula which was known for hyperelliptic curves and recently proved for superelliptic curves \cite{book} is not known to be true in general for algebraic curves.

In preparing this paper we used as a blueprint a talk of the third author from 2012.  The  excellent survey of Griffiths \cite{gri} and the publication of \cite{book} were timely in helping us understand some of the classical works from the more modern viewpoint.  We have tried to give somewhat a historical perspective on theta functions including references to the work of Abel, Jacobi, Weierstrass, Riemann, and the modern viewpoint of Mumford and others.

The story starts with Legendre, who used elliptic functions for problems such as the movement of a simple pendulum and the deflection of a thin elastic bar.
Legendre spent more than forty years of his life working on elliptic functions, including the classification of elliptic integrals and published the following  books: 
\textit{Elliptic integrals} (1786), 
\textit{Elliptic transcendents} (1792),  and \textit{Elliptic functions} (1811-1816)  which appeared in 3 volumes.
Despite forty years of dedication to elliptic functions, Legendre's work went essentially unnoticed by his contemporaries until   Abel and Jacobi's work on the subject. 

In 1825, the Norwegian government funded Abel on a scholarly visit to France and Germany. Abel then traveled to Paris, where he gave an important paper revealing the double periodicity of the
elliptic functions.

Jacobi (1829) wrote the classic treatise on elliptic functions, of great importance in mathematical physics, because of the need to integrate second order kinetic energy equations. The motion equations in rotational form are integrable only for the three cases of the pendulum, the symmetric top in a gravitational field, and a freely spinning body, wherein solutions are in terms of elliptic functions. Jacobi was also the first mathematician to apply elliptic functions to number theory, for example, proving the polygonal number theorem of Pierre de Fermat.

In developments of the theory of elliptic functions, modern authors mostly follow Karl Weierstrass. The notations of Weierstrass  elliptic functions based on his p-function are convenient, and any
elliptic function can be expressed in terms of these. The elliptic functions introduced by Carl Jacobi, and the auxiliary theta functions (not doubly-periodic), are more complex but important both for the
history and for general theory.

Riemann (1826-1866) developed the general theory of theta functions generalizing the Jacobi's theta functions.
Riemann's dissertation, completed under Gauss's supervision in 1851, was on the foundations of complex analysis. It introduced several ideas of fundamental importance, such as the definitions of conformal mapping and simple connectivity. These are necessary for one of his main results, the Riemann mapping theorem: any simply connected domain of the complex plane having at least two boundary points can be conformally mapped onto the unit disk. Riemann also introduced the Laurent series expansion for functions having poles and branch points.

In Section \ref{sect-2},  we give a short historical view of what led to the definition of the theta functions. We  describe  elliptic integrals, Abelian integrals, and the Abel's theorem.  A wonderful source on such topics is Baker's book \cite{baker} among more modern viewpoints of the area. A nice article discussing some of the historical aspects and Abel's contribution is Griffiths paper \cite{gri} on the bicentenary of Abel's birthday. 
We continue with a overly simplified version of   the Jacobi inversion problem which lead to elliptic functions and the Jacobi theta functions.

In section \ref{sect-3} we define the Riemann-theta functions, periods, characteristics, and G\"opel systems. For each G\"opel system we find identities among theta-nulls as in \cite{Sh-1} and \cite{Sh-2}. In section \ref{sect-4}, the focus shifts to the hyperelliptic curves.  In such case determining identities among thetanulls can be worked out explicitly due to Frobenius' theta formula (cf. Lem.~\ref{Frob}) and the Thomae's formula (cf. Lem.~\ref{Thomae}).  Both such important results are described in detail.  

Furthermore, we compute explicitly such identities for genus 2 curves as in \cite{Sh-2}, give the Picard's lemma  (cf. Lem.~\ref{pic}) for genus 2 curves, and prove an analogues of the Picard's lemma where the branch points are expressed only in terms of the fundamental thetanulls $\th_1, \th_2, \th_3, \th_4$ (cf. Lem.\ref{possibleCurve}).  Moreover, we give algebraic relations among the fundamental theta functions when the curve has at least two non-hyperelliptic involutions (cf. Thm.~\ref{V4}).  
We perform similar computations for genus 3 hyperelliptic curves even though computations are longer and less explicit.  It would be interesting to see if such computations become easier in the view of the absolute invariants of genus 3 hyperelliptic curves as in \cite{hyp-3}.

In section \ref{sect-5}, we shift our attention to superelliptic curves.  The goal of determining identities among thetanulls is not easily achieved because the Thomae's formula is not as simple as in the case of hyperelliptic case.  There has been a lot of activity in the last decade on proving an analogue of the Thomae's formula for superelliptic curves; see \cite{Enolski}, \cite{NK}, and others and summarized in \cite{book}. We give a version of the Thomae's formula in Thm.~\ref{CyclicCurvesThm}.

In section \ref{sect-6},   we develop an algorithm to determine relations among theta functions of a cyclic curve $\X$ with automorphism group $\Aut(\X)$.

There are many unanswered questions when it comes to theta functions for curves of higher genus.  However, theoretically we now know how to determine the fundamental theta functions for superelliptic curves of any genus. 

\medskip

\noindent \textbf{Notation:} Throughout this paper $k$ denotes an algebraically closed field of characteristic zero, $g$ an integer $\geq 2$, and $\X_g$ a hyperelliptic curve of genus $g$ defined over $k$. A cyclic group of order $n$ will be denoted by $C_n$, unless otherwise stated. 

\section{Abelian integrals, some historical remarks}\label{sect-2}
An \textit{algebraic function} $y(x)$ is a function which satisfies some equation 
\[ f(x, y(x)) =0, \]
where $f(x, y) \in \C[x, y]$ is an irreducible polynomial. In the beginning of XIX century it was a lot of interest for integrals of algebraic functions, as mentioned in the introduction.  Below we give a brief description of such integrals which lead to theorems of Abel and Jacobi. 

Recall from Calculus that $\int F(x) \, dx$, 
for $F(x) \in \C (x)$,  can be solved via the partial fractions method by expressing this as a sum of rational functions in $x$ or logarithms of $x$.   
Also, the integral  
\[ \int F(x, y) \, dx, \] 
where $F \in \C (x, y)$ and $x, y \in \C(t)$, can be easily solved   by replacing for $x=x(t)$ and $y=y(t)$ this reduces to the previous case.   Similarly, we can deal with the case   
\[ \int F \left(x,  \sqrt{ax^2+bx+c} \right) \, dx. \]
Indeed, let $y=\sqrt{ax^2+bx+c}$.  Then, $y^2=ax^2+bx+c$ is the equation of a conic.  As such it can be parametrized as $x= x(t)$,  $y=y(t)$ 
and again reduces to the previous case.  

However, the integral   
\[ \int F \left(x,  \sqrt{ax^3+bx^2+cx+d} \right) \, dx \]  
can not be solved this way because  
\[ y^2= ax^3+bx^2+cx+d\] 
is not a genus 0 curve, and therefore can not be parametrized.    Such integrals are called \textbf{elliptic integrals}.  To solve them one needs to understand the concept of  \textbf{elliptic functions} which will be developed later.  It can be easily shown that these integrals can be transformed to the form 
\[ \int \frac {p(x)} {\sqrt{q(x)} } dx \]
where $p(x), q(x)$ are polynomials such that $\deg q = 3, 4$ and $q(x)$ is separable.  The term \textit{elliptic} comes from the fact that such integrals come up in the computation of the length of an ellipse.

\begin{exe}
Let an ellipse be given by 
\[ \frac {x^2} {a^2} + \frac {y^2} {b^2} =1, \quad a> b,\]
and denoted by $k^2 = \frac {a^2-b^2} {a^2}$.  We denote by  $t = \arcsin \frac x a$. 
Prove that the arc length of the ellipse is given by 
\[ L = a \int_0^{2\pi} \frac {1-k^2x^2} {\sqrt{(1-t^2) (1-k^2 t^2}} \, dt \]
\end{exe}

It is worth reminding our reader that elliptic integrals are the first when we can not solve them via the elementary calculus. In other words, they can not be expressed as a sum of rational and logarithmic functions.  Instead, we need other transcendental functions, namely the elliptic functions. 


A natural generalization of the elliptic integrals are the \textbf{hyperelliptic integrals} which are of the form 
\[ \int \frac {p(x)} {\sqrt{q(x)} } dx \]
where $p(x), q(x)$ are polynomials such that $\deg q \geq 5$ and $q(x)$ is separable.

Naturally, the square root above can be assumed to be a n-th root.  We will call such integrals \textbf{superelliptic integrals}.  Hence, a superelliptic integral is of the form \[ \int \frac {p(x)} {\sqrt[n]{q(x)} } dx \]
where $n \geq 3$, $p(x), q(x)$ are polynomials such that $\deg q \geq 5$ and $q(x)$ is separable. 

What about the general case when 
\[ \int R(x, y) \, dx, \]  where $R \in \C (x, y)$ and $y$ is an algebraic function of $x$ given by some equation $F(x, y) =0$, 
for $F(x, y)\in \C [x, y]$?   An integral of this type is called an \textbf{Abelian integral}.  

In the next few sections we will describe how the theory of Abelian integrals led to some fundamental results in mathematics and its role in developing of algebraic geometry.  

\subsection{Abel's theorem}

There are several version of what is called the Abel's theorem in the literature. For original versions of what Abel actually stated and proved one can check the classic books \cite{baker} and \cite{Clebsch}.  For modern interpretations of Abel's theorem and its historical perspectives there are the following wonderful references \cite{gri},  \cite{gri-2} and \cite{kleiman}. In this short notes we will try to stay as close as possible to the original version of Abel.  

Let $y$ be an algebraic function of $x$ defined by an equation of the form 
\[ f(x, y) = y^n + A_1, y^{n-1} + \cdots A_n=0,\]
where $A_0, \dots , A_n \in \C (x)$. Let $R (x, y) \in \C(x, y)$. 

\begin{thm}[Abel]
The sum \[ \int_{(a_1, b_1)}^{(x_1, y_1)} R(x, y) + \cdots + \int_{(a_m, b_m)}^{(x_m, y_m)} R(x, y) \]
for arbitrary $a_i, b_i$, is expressible as a sum of rational functions of $(x_1, y_1)$, $\dots$, $(x_m, y_m)$ and logarithms of such rational functions with the addition of 
\[ - \int^{(z_1, s_1)} R(x, y)  - \cdots - \int^{(z_k, s_k)} R(x, y)\]
where $z_i, s_i$ are determined by $x_i, y_i$ as the roots of an algebraic equation whose coefficients are rational coefficients of $x_1, y_1, \dots , x_m, y_m$ and $s_1, \dots , s_k$ are the corresponding values of $y$, for which any $s_i$ is determined as a rational function of $z_i$ and   $x_1, y_1, \dots , x_m, y_m$.
Moreover, the number $k$ does not depend on $m$, $R(x, y)$, or the values $(x_i, y_i)$, but only on the equation 
\[ f(x, y) =0.\]
\end{thm}
For more details of this version of Abel's theorem and its proof see \cite[pg. 207-235]{baker}.
A modern version of the Abel's theorem, which is found in most textbooks says that the Abel-Jacobi's map is injective; see Thm.~\ref{thm-abel} for details. A nice discussion from the modern point of view is \cite{kleiman}. 

\subsection{Jacobi inversion problem}

The new idea of Jacobi was to consider integrals $\int_{c}^{w} R(x, y)$ as variables and to try to determine $w$ in terms of such variables. This idea led to the fundamental concept of theta functions, which will be formally defined in the next section.

First, consider the Abelian integrals  \[  \int_{c_i}^{w_i} R(x, y) = z_i\]  for $i=1, \dots g$.  Consider 
\[  z_i : =\int_{c_i}^{w_i} R(x, y) \] 
as variables and express $w_i$ as functions of $z_i$,
\[ w_i = f(z_i).\]
This is known as the \textit{Jacobi inversion problem}.  


\begin{exa}[Elliptic integrals]   Let be given the integral (i.e. $g=1$)
\[ \int_0^{w_1} \frac {dt}   {\sqrt{(1-t^2)(1-k^2t^2) } } = z_1\]
Then \[w_1= sn (z_1)= sn (u; k) = \frac  {\th_3(0) \th_1 (v)} {\th_2 (0) \th_0 (v)}, \]   where $u=v\, \pi \, \th_3^2 (0)$ and $\th_0, \th_1, \th_2, \th_3$ are the Jacobi theta functions; see \cite{baker} for details. 
\end{exa}
%

It was exactly the above case that was the motivation of Jacobi to introduce the theta functions.  With these functions he expressed his functions $\sn u$, $\cn u$, and $\dn u$ as fractions having the same denominators, with zeroes of this denominator being the common poles of $\sn u$, $\cn u$, and $\dn u$. 

For  $g=2$, G\"opel found similar functions, building on work of Hermite.  We will say more about this case in the coming sections. G\"opel and later Rosenhain notice that integrals of the first kind, which exist for $g=2$ become elliptic integrals of the first and third kind, when two branch points of the curve of $g=2$ coincide.  This case corresponds to the degenerate cases of the $\L_n$ spaces as described in \cite{phd} and later in \cite{deg3}. Both G\"opel and Rosenhain in developing theta functions for genus $g=2$ were motivated by the Jacobi inversion problem.  Weierstrass considered functions which are quotients of theta functions for the  hyperelliptic curves, even though it seems as he never used the term "theta functions".  

In their generality, theta functions were developed by Riemann for any $g \geq 2$.  It is Riemann's approach that is found in most modern books and that we will briefly describe in the next section. Most known references for what comes next can be found in \cite{Igusa, Mu1, Mu2, Mu3}.

\section{Riemann's theta functions}\label{sect-3}
In this section we define the Riemann-theta functions, theta characteristics, and theta-nulls which will be the main focus for the rest of the paper.  Most of the material in this section is taken from \cite{Sh-2}. 

\subsection{Introduction to theta functions of curves}
Let $\X$ be an irreducible, smooth, projective curve of genus $g \geq 2$ defined over the complex field $\C.$ We denote the moduli space of genus $g$ by $\M_g$ and  the hyperelliptic locus in $\M_g$ by $\H_g.$ It is well known that  $\dim \M_g = 3g-3$ and $\H_g$ is a $(2g-1)$ dimensional subvariety of $\M_g.$  

Choose a symplectic homology basis for $\X$, say
\[ \{ A_1, \dots, A_g, B_1, \dots , B_g\}\]
such that the intersection products $A_i \cdot A_j = B_i \cdot B_j =0$ and $A_i \cdot B_j= \d_{i j}$.
We choose a basis $\{ w_i\}$ for the space of holomorphic 1-forms such that $\int_{A_i} w_j = \d_{i j},$ where $\d_{i j}$ is the Kronecker delta. The matrix $\O= \left[ \int_{B_i} w_j \right] $ is  the \textbf{period matrix} of $\X$.
The columns of the matrix $\left[ I \ | \O \right]$ form a lattice $L$ in  $\C^g$ and the Jacobian  of $\X$ is $\J(\X)= \C^g/ L$. 


Fix a point $p_0 \in \X$.  Then, the Abel-Jacobi map is defined as follows
\[
\begin{split}
\mu_p : \X & \to \J (\X) \\
p & \to \left( \int_{p_0}^p w_1, \dots ,   \int_{p_0}^p w_g \right) \mod L \\
\end{split}
\]
The Abel-Jacobi map can be extended to divisors of $\X$ the natural way, for example for a divisor $D= \sum_i n_i P_i$ we defined %
\[ \mu (D) = \sum_i n_i \mu (P_i). \]
The following two theorems are part of the folklore on the subject and their proofs can be found in all classical textbooks. 
\begin{thm}[Abel]\label{thm-abel}
The Abel-Jacobi map is injective.
\end{thm}

\begin{thm}[Jacobi]\label{thm-jacobi}
The Abel-Jacobi map is surjective
\end{thm}

We continue with our goal of defining  theta functions and theta characteristics.   Let 
\[ \HS_g =\{\t : \t \,\, \textit{is symmetric}\,\, g \times g \, \textit{matrix with positive definite imaginary part} \}\]
be the \textbf{Siegel upper-half space}. Then $\O \in \HS_g$. The group of all $2g \times 2g$ matrices $M \in GL_{2g}(\Z)$ satisfying
\[M^t J M = J  \,\,\,\,\,\,\,\, \textit{with} \, \,\,\,\,\,\, J = \begin{pmatrix}  0 & I_g \\ -I_g & 0 \end{pmatrix} \]
is called the \textbf{symplectic group} and denoted  by $Sp_{2g}(\Z)$.
Let $M = \begin{pmatrix} R & S \\ T & U \end{pmatrix} \in Sp_{2g}(\Z) $ and $\t \in \HS_g$ where $R,$ $S,$ $T$ and $U$ are $g \times g$ matrices. 
$Sp_{2g}(\Z) $ acts transitively on $\HS_g$ as

\[ M(\t) = (R \t + S)(T \t + U)^{-1}. \]

\noindent Here, the multiplications are matrix multiplications. There is an injection
\[ \M_g \embd \HS_g/ Sp_{2g}(\Z) =: \A_g, \] 
where each curve $C$ (up to isomorphism) goes to its Jacobian in $\A_g.$

If $\ell$ is a positive integer, the principal congruence group of degree $g$ and of level $\ell$ is defined as a subgroup of $Sp_{2g}(\Z)$ by the condition $M \equiv I_{2g} \mod \ell.$ We shall denote this group by $Sp_{2g}(\Z)(\ell)$.

For any $z \in \C^g$ and $\t \in \HS_g$ the \textbf{Riemann's theta function} is defined as
\[ \T (z , \t) = \sum_{u\in \Z^g} e^{\pi i ( u^t \t u + 2 u^t z )  }\]
where $u$ and $z$ are $g$-dimensional column vectors and the products involved in the formula are matrix products. The fact that the imaginary part of $\t$ is positive makes the series absolutely convergent over every compact subset of $\C^g \times \HS_g$.

The theta function is holomorphic on $\C^g\times \HS_g$ and has quasi periodic properties,
\[ \T(z+u,\tau)=\T(z,\tau)\quad \textit{and}\quad \T(z+u\tau,\tau)=e^{-\pi i( u^t \tau u+2z^t u )}\cdot  \T(z,\tau), \]
where $u\in \Z^g$; see \cite{Mu1} for details.  The locus 
\[ \Theta: = \{ z \in \C^g/L : \T(z, \O)=0 \}\]
is called the \textbf{theta divisor} of $\X$.
Any point $e \in \J (\X)$ can be uniquely written  as $e = (b,a) \begin{pmatrix} 1_g \\ \O \end{pmatrix}$ where $a,b \in \R^g$ are the characteristics of $e.$
We shall use  the notation $[e]$ for the characteristic of $e$ where $[e] = \ch{a}{b}.$ For any $a, b \in \Q^g$, the theta function with rational characteristics is defined as a translate of Riemann's theta function multiplied by an exponential factor

\begin{equation} \label{ThetaFunctionWithCharac} \T  \ch{a}{b} (z , \t) = e^{\pi i( a^t \t a + 2 a^t(z+b))} \T(z+\t a+b ,\t).\end{equation}
\noindent By writing out Eq.~\eqref{ThetaFunctionWithCharac}, we have
\[ \T  \ch{a}{b} (z , \t) = \sum_{u\in \Z^g} e^{\pi i ( (u+a)^t \t (u+a) + 2 (u+a)^t (z+b) )  }. \]
The Riemann's theta function is $\T \ch{0}{0}.$ The theta function with rational characteristics has  the following properties:
\begin{equation}\label{periodicproperty}
\begin{split}
& \T \ch{a+n} {b+m} (z,\t) = e^{2\pi i a^t m}\T \ch {a} {b} (z,\t),\\
&\T \ch{a} {b} (z+m,\t) = e^{2\pi i a^t m}\T \ch {a} {b} (z,\t),\\
&\T \ch{a} {b} (z+\t m,\t) = e^{\pi i (-2b^t m -m^t \t m - 2m^t z)}\T \ch {a} {b} (z,\t)\\
\end{split}
\end{equation}
where $n,m \in \Z^n.$ All of these properties are immediately verified by writing them out.

A scalar obtained by evaluating a theta function with characteristic at $z=0$ is called a \emph{theta constant} or \emph{theta-nulls}. When the entries of column vectors $a$ and $b$ are from the set $\{ 0,\frac{1}{2}\}$, then the characteristics $ \ch {a}{b} $ are called the \emph{half-integer characteristics}. The corresponding theta functions with rational characteristics are called \emph{theta characteristics}.

Points of order $n$ on $\J(\X)$ are called the $\frac 1 n$-\textbf{periods}. Any point $p$ of $\J(\X)$ can be written as $p = \t \,a + b. $ If $\ch{a}{b}$ is a $\frac 1 n$-period, then $a,b \in (\frac{1}{n}\Z /\Z)^{g}.$ The $\frac 1 n$-period $p$ can be associated with an element of $H_1(\X,\Z / n\Z)$ as follows: 

Let $a = (a_1,\cdots,a_g)^t,$ and $b = (b_1,\cdots,b_g)^t.$ Then
\[
\begin{split}
p  & = \t a + b \\
           & = \left(      \sum a_i \int_{B_i} \om_1, \cdots , \sum a_i \int_{B_i} \om_g \right)^t     +       
           \left(b_1 \int_{A_1} \om_1, \cdots  ,   b_g       \int_{A_g} \om_g \right) \\
           & = \left(\sum (a_i \int_{B_i} \om_1 + b_i\int_{A_i} \om_1 \right) , \cdots , \sum \left(a_i \int_{B_i} \om_g + b_i\int_{A_i} \om_g) \right)^t\\
           & = \left( \int_C \om_1, \cdots, \int_C \om_g \right)^t
\end{split}
\]
where $C = \sum a_i B_i + b_i A_i. $ We identify the point $p$ with the cycle $\bar{C} \in H_1(\X,\Z / n\Z)$ where $\bar{C} =\sum \bar{a_i} B_i + \bar{b_i} A_i,$  $\bar{a_i} = n a_i$ and $\bar{b_i} = n b_i$ for all $i$;  see  \cite{Accola} for more details. 

\subsubsection{Half-Integer Characteristics and the G\"opel Group} 
In this section we study groups of half-integer characteristics. Any half-integer characteristic $\m \in\frac{1}{2}\Z^{2g}/\Z^{2g}$ is given by
\[
\m = \frac{1}{2}m = \frac{1}{2}
\begin{pmatrix} m_1 & m_2 &  \cdots &  m_g \\ m_1^{\prime} & m_2^{\prime} & \cdots & m_g^{\prime}  \end{pmatrix},
\]
where $m_i, m_i^{\prime} \in \Z.$ For $\m = \ch{m ^\prime}{m^{\prime \prime}} \in \frac{1}{2}\Z^{2g}/\Z^{2g},$ we define $e_*(\m) = (-1)^{4 (m^\prime)^t m^{\prime \prime}}.$ We say that $\m$ is an \emph{even} (resp. \emph{odd}) characteristic if $e_*(\m) = 1$ (resp. $e_*(\m) = -1$). For any curve of genus $g$, there are $2^{g-1}(2^g+1)$ (resp., $2^{g-1}(2^g-1)$ ) even theta functions (resp., odd theta functions). Let $\an$ be another half-integer characteristic. We define
\[
 \m \, \an = \frac{1}{2} \begin{pmatrix} t_1 & t_2 &  \cdots &  t_g \\ t_1^{\prime} & t_2^{\prime} & \cdots &
t_g^{\prime}
\end{pmatrix}
\]
where $t_i \equiv (m_i\, + a_i)  \mod 2$ and $t_i^{\prime} \equiv (m_i^{\prime}\, + a_i^{\prime} ) \mod 2.$

For the rest of this paper we only consider  characteristics $\frac{1}{2}q$ in which each of the elements
$q_i,q_i^{\prime}$ is either 0 or 1. We use the following abbreviations:
\[
\begin{split}
&|\m| = \sum_{i=1}^g m_i m_i^{\prime},  \quad \quad \quad \quad \quad \quad \quad \quad \quad
|\m, \an| = \sum_{i=1}^g (m_i^{\prime} a_i - m_i a_i^{\prime}), \\
& |\m, \an, \bn| = |\an, \bn| + |\bn, \m| + |\m, \an|, \quad \quad {\m\choose \an} = e^{\pi i \sum_{j=1}^g m_j
a_j^{\prime}}.
\end{split}
\]
\indent The set of all half-integer characteristics forms a group $\G$ which has $2^{2g}$ elements. We say that two half integer characteristics $\m$ and $\an$ are \emph{syzygetic} (resp., \emph{azygetic}) if $|\m, \an| \equiv 0 \mod 2$ (resp., $|\m, \an| \equiv 1 \mod 2$) and three half-integer characteristics $\m, \an$, and $\bn$ are syzygetic if
$|\m, \an, \bn| \equiv 0 \mod 2$.

A \emph{G\"opel group} $G$ is a group of $2^r$ half-integer characteristics where $r \leq g$ such that every two characteristics are syzygetic. The elements of the group $G$ are formed by the sums of $r$ fundamental characteristics; see \cite[pg. 489]{baker} for details. Obviously, a G\"opel group of order $2^r$ is isomorphic to $C^r_2$. The proof of
the following lemma can be found on   \cite[pg.  490]{baker}.
\begin{lem}
The number of different G\"opel groups which have $2^r$ characteristics is
\[ \frac{(2^{2g}-1)(2^{2g-2}-1)\cdots(2^{2g-2r+2}-1)}{(2^r-1)(2^{r-1}-1)\cdots(2-1)}. \]
\end{lem}
If $G$ is a G\"opel group with $2^r$ elements, it has $2^{2g-r}$ cosets. The cosets are called \emph{G\"opel systems}
and are denoted by $\an G$, $\an \in \G$. Any three characteristics of a G\"opel system are syzygetic. We can find a
set of characteristics called a basis of the G\"opel system which derives all its $2^r$ characteristics by taking only
combinations of any odd number of characteristics of the basis.
\begin{lem}
Let $g \geq 1$ be a fixed integer, $r$ be as defined above and $\sigma = g-r.$ Then there are
$2^{\sigma-1}(2^\sigma+1)$ G\"opel systems which only consist of even characteristics and there are
$2^{\sigma-1}(2^\sigma-1)$ G\"opel systems which consist of odd characteristics. The other $2^{2\sigma}(2^r-1)$ G\"opel
systems consist of as many odd characteristics as even characteristics.
\end{lem}
\proof The proof can be found on \cite[pg. 492]{baker}. \qed
\begin{cor}\label{numb_systems}
When $r=g,$ we have only one (resp., 0) G\"opel system which consists of even (resp., odd) characteristics.
\end{cor}
Let us consider $s=2^{2\sigma}$ G\"opel systems which have  distinct characters. Let us denote them by
\[\an_1 G,\an_2 G,\cdots,\an_s G.\] We have the following lemma.
\begin{lem}
It is possible to choose $2\sigma+1$ characteristics from $\an_1, \an_2,\cdots, \an_s,$  say $\bar{\an}_1,$
$\bar{\an}_2,$ $\cdots,$ $\bar{\an}_{2\sigma+1}$, such that every three of them are azygetic and all have the same
character. The above $2\sigma+1$ fundamental characteristics are even (resp., odd) if $\sigma \equiv 1,0 \mod 4$
(resp.,$\equiv 2,3 \mod 4$).
\end{lem}
\noindent The proof of the following lemma can be found on \cite[pg. 511]{baker}.

\begin{lem}
For any half-integer characteristics $\an$ and $\hn,$ we have the following:
\begin{equation}\label {Bakereq1}
\T^2[\an](z_1,\t) \T^2[\an \hn](z_2,\t) = \frac{1}{2^{g}} \sum_\en  e^{\pi i |\an \en|} { \hn \choose \an \en}
\T^2[\en](z_1,\t)\T^2[\en \hn](z_2,\t).
\end{equation}
\end{lem}

We can use this relation to get identities among half-integer thetanulls. Here $\en$ can be any half-integer
characteristic. We know that we have $2^{g-1}(2^g+1)$ even characteristics. As the genus increases, we have multiple
choices for $\en.$ In the following, we explain how we reduce the number of possibilities for $\en$ and how to get
identities among thetanulls.

First we replace $\en$ by $\en \hn$ and $z_1=z_2= 0$ in Eq.~\eqref{Bakereq1}. Eq.~\eqref{Bakereq1} can then be written
as follows:
\begin{equation}\label {Bakereq2}
\T^2[\an] \T^2[\an \hn] = 2^{-g} \sum_\en  e^{\pi i |\an \en \hn|} { \hn \choose \an \en \hn} \T^2[\en] \T^2[\en \hn].
\end{equation}
We have $e^{\pi i |\an \en \hn|}{ \hn \choose \an \en \hn} = e^{\pi i |\an \en|}{ \hn \choose \an \en} e^{\pi i |\an
\en, \hn|}.$ Next we put $z_1=z_2= 0$ in Eq.~\eqref{Bakereq1} and add it to Eq.~\eqref{Bakereq2} and get the following
identity:
\begin{equation}\label {Bakereq3}
2\T^2[\an] \T^2[\an \hn] = 2^{-g} \sum_\en  e^{\pi i |\an \en|} (1 + e^{\pi i|\an \en, \hn|}) \T^2[\en] \T^2[\en \hn].
\end{equation}
If $|\an \en, \hn| \equiv  1 \mod 2$, the corresponding terms in the summation vanish. Otherwise $1 + e^{\pi i|\an \en,
\hn|} = 2.$ In this case, if either $\en$ is odd or $\en \hn$ is odd, the corresponding terms in the summation vanish
again. Therefore, we need $|\an \en, \hn| \equiv 0 \mod 2$ and $|\en| \equiv |\en \hn| \equiv 0 \mod 2,$ in order to
get nonzero terms in the summation. If $\en^*$ satisfies $|\en^*| \equiv |\en^* \hn^*| \equiv 0 \mod 2$ for some
$\hn^*,$ then $\en^*\hn^*$ is also a candidate for the left hand side of the summation. Only one of such two values
$\en^*$ and $\en^* \hn^*$ is taken. As  a result, we have the following identity among thetanulls
\begin{equation}\label {eq1}
\T^2[\an] \T^2[\an \hn] = \frac{1}{2^{g-1}} \sum_\en  e^{\pi i |\an \en|} { \hn \choose \an \en} \T^2[\en]\T^2[\en
\hn],
\end{equation}
where $\an, \hn$ are any characteristics and $\en$ is a characteristics such that $|\an \en, \hn| \equiv 0 \mod 2,$
$|\en| \equiv |\en \hn| \equiv 0 \mod 2$ and $\en \neq \en \hn.$

By starting from the Eq.~\eqref{Bakereq1} with $z_1 = z_2$ and following a similar argument to the one above, we can
derive the identity,
\begin{equation}\label{eq2}
\T^4[\an] + e^{\pi i |\an, \hn|} \T^4[\an \hn] = \frac{1}{2^{g-1}} \sum_\en  e^{\pi i |\an \en|} \{ \T^4[\en] + e^{ \pi
i |\an, \hn|} \T^4[\en \hn]\}
\end{equation}
where $\an, \hn$ are any characteristics and $\en$ is a characteristic such that $|\hn| + |\en, \hn| \equiv 0 \mod 2,$
$|\en| \equiv |\en \hn| \equiv 0 \mod 2$ and $\en \neq \en \hn.$

\begin{rem}
$|\an \en ,\hn| \equiv 0 \mod 2$ and $|\en \hn| \equiv |\en| \equiv 0 \mod 2$ implies $|\an, \hn| + |\hn| \equiv 0 \mod
2.$
\end{rem}
We use Eq.~\eqref{eq1} and Eq.~\eqref{eq2} to get identities among theta-nulls.
%


\section{Hyperelliptic curves and their theta functions}\label{sect-4}
A hyperelliptic curve $\X,$ defined over $\C,$ is a cover of order two of the projective line $\P^1.$  
 Let $\X \longrightarrow \P^1$ be the degree 2 hyperelliptic projection.
We can assume that $\infty$ is a branch point.

Let  $B := \{\a_1,\a_2, \cdots ,\a_{2g+1} \}$
be the set of other branch points and let $S = \{1,2, \cdots, 2g+1\}$ be the index set of $B$ and $\e : S \longrightarrow
\frac{1}{2}\Z^{2g}/\Z^{2g}$  be a map defined as follows:
%
\[
\begin{split}
\e(2i-1) & = \begin{bmatrix}
              0 & \cdots & 0 & \frac{1}{2} & 0 & \cdots & 0\\
              \frac{1}{2} & \cdots & \frac{1}{2} & 0 & 0 & \cdots & 0\\
            \end{bmatrix}, \\
 \e(2i) & =\begin{bmatrix}
              0 & \cdots & 0 & \frac{1}{2} & 0 & \cdots & 0\\
              \frac{1}{2} & \cdots & \frac{1}{2} & \frac{1}{2} & 0 & \cdots & 0\\
            \end{bmatrix}
\end{split}
\]
%
where the nonzero element of the first row appears in $i^{th}$ column.
We define  $\e(\infty) $ to be $
\begin{bmatrix}
              0 & \cdots & 0 & 0\\
              0 & \cdots & 0 & 0\\
            \end{bmatrix}$.
For any $T \subset B $, we define the half-integer characteristic as

\[ \e_T = \sum_{a_k \in T } \e(k) .\]
Let $T^c$ denote the complement of $T$ in $B.$ Note that $\e_B \in \Z^{2g}.$ If we view $\e_T$ as an element of
$\frac{1}{2}\Z^{2g}/\Z^{2g}$ then $\e_T= \e_{T^c}.$ Let $\triangle$ denote the symmetric difference of sets, that is $T
\triangle R = (T \cup R) - (T \cap R).$ It can be shown that the set of subsets of $B$ is a group under $\triangle.$ We
have the following group isomorphism:
\[ \{T \subset B\,  |\, \#T \equiv g+1 \mod 2\} / T \sim T^c \cong \frac{1}{2}\Z^{2g}/\Z^{2g}.\]
For $\gamma = \ch{\gamma ^\prime}{\gamma^{\prime \prime}} \in \frac{1}{2}\Z^{2g}/\Z^{2g}$, we have
\begin{equation}\label{parityIdentity}
 \T [\gamma] (-z , \t) = e_* (\gamma) \T [\gamma] (z , \t).\end{equation}
It is known that for hyperelliptic curves, $2^{g-1}(2^g+1) - {2g+1 \choose g}$ of the even thetanulls are zero.
The following theorem provides a condition for the characteristics in which theta characteristics become zero. The
proof of the theorem can be found  in \cite{Mu2}.
\begin{thm}\label{vanishingProperty}
Let $\X$ be a hyperelliptic curve, with a set $B$ of branch points. Let $S$ be the index set as above and $U $ be the
set of all odd values of $S$. Then for all $T \subset S$ with even cardinality, we have $ \T[\e_T] = 0$  if and only if
$\#(T \triangle U) \neq g+1$, where $\T[\e_T]$ is the theta constant corresponding to the characteristics $\e_T$.
\end{thm}
When the characteristic $\gamma$ is odd, $e_* (\gamma)=1.$ Then from Eq.~\eqref{parityIdentity} all odd thetanulls
are zero. There is a formula which satisfies half-integer theta characteristics for hyperelliptic curves called
\textit{Frobenius' theta formula}.

\begin{lem}[Frobenius]\label{Frob}
For all $z_i \in \C^g$, $1\leq i \leq 4$ such that $z_1 + z_2 + z_3 + z_4 = 0$ and for all $b_i \in \Q^{2g}$, $1\leq i
\leq 4$ such that $b_1 + b_2 + b_3 + b_4 = 0$, we have
\[ \sum_{j \in S \cup \{\infty\}} \epsilon_U(j) \prod_{i =1}^4 \T[b_i+\e(j)](z_i) = 0, \]
where for any $A \subset B$,
\[
\epsilon_A(k) =
          \begin{cases}
           1 & \textit {if $k \in A$}, \\
           -1 & \textit {otherwise}.
          \end{cases}
\]
\end{lem}
\proof See \cite[pg.107]{Mu1}. \qed

A relationship between thetanulls and the branch points of the hyperelliptic curve is given by Thomae's formula.
\begin{lem}[Thomae]\label{Thomae}
%
For all sets of branch points $B=\{\a_1,\a_2, \cdots ,\a_{2g+1} \},$ there is a constant $A$ such that for all $T\subset B,$ $\# T$ is even, \\
\[\T[\eta_T](0;\t)^4 =(-1)^{\#T \cap U} A \prod_{i<j \atop i,j \in T \triangle U} (\a_{i} - \a_{j}) \prod_{i<j \atop i,j \notin T \triangle U} (\a_{i} -
\a_{j})  \] \\
\noindent where $\eta_T$ is a non singular even half-integer characteristic corresponding to the subset $T$ of branch
points.
\end{lem}
See \cite[pg. 128]{Mu1} for the description of $A$ and \cite[pg.  120]{Mu1} for the proof. Using Thomae's formula and Frobenius' theta identities we express the branch points of the hyperelliptic curves in terms of even thetanulls. 
In \cite{Sh-1} and \cite{Sh-2} it is shown how such relations are computed for genus $g=2, 3$.

\subsection{Superelliptic  curves and their  theta functions}\label{sect-5}

Generalizing the theory of theta functions of hyperelliptic curves to all cyclic covers of the projective line has been the focus of research of the last few decades.  The main efforts have been on generalizing the Thomae's formula to such curves.  In the literature of Rimann surfaces such curves are called for historical reasons the $\Z_n$ curves. For a summary of some of the results on the Thomae's formula for $\Z_n$ curves and especially the relations of thetanulls for such curves with extra automorphisms the reader can check \cite{Sh-1}, \cite{Kop}, \cite{W}, \cite{Sh-2}. Especially in \cite{W} and  \cite{Sh-2} are summarized the known results up to that time \cite{br}, \cite{Enolski}, \cite{NK}, cite{SHI}.

As a more recent development came out a book in this topic \cite{book}. Obviously it would be a difficult task for us to sumarize all the results of \cite{book} in this short section.  A condensed account fo that philosophy and some computational results can be found in Wijesiri's thesis \cite{W}. 

\section{Vanishing of theta nulls for genus 3 curves with automorphisms}
In this section we focus on the genus 3 curves with the goal of describing the loci $\M_3 (G, \bC)$ in terms of theta functions for each possible group $G$ and signature $\bC$. For the rest of this paper $\X$ denotes a genus 3 algebraic curve defined over $\C$.

A covering $f :  X  \to Y$ of algebraic varieties is called a \textbf{maximal covering}   if it does not factor over a nontrivial isogeny. A map of algebraic curves $f: X \to Y$ induces maps between their Jacobians 
$f^*: \J (Y) \to \J (X)$ and $f_*: \J (X) \to \J (Y)$. When $f$ is maximal then $f^*$ is injective and $\ker (f_*)$ is connected, see \cite{Se} (p. 158) for details.

Our strategy is to find an appropriate element $\s \in \Aut (\X)$ and study the cover $\pi: \X \to \X/\<\s\>$. We will denote the quotient space $\X/\<\s\>$ by $\X^\s$. Studying the Jacobian $\J (\X^\s)$ and using the induced map $\pi^\ast : \J (\X^\s) \to \J (X)$ we would like to say something about the $\frac 1 n $ --periods of $\J (\X)$.

Next we recall a classical result on half-periods; See Krazer \cite[pg. 294, XXXII Satz]{Kr} for the proof.

\begin{prop} Let $G \leq \J (\X)$ and $G$ is generated by distinct half-periods $G:=\< \s_1, \dots , \s_r\> $. Then $\X$ has a basis $\{\a_1, \dots , \a_m, \b_1, \dots \b_{2n}\}$ with $m + 2n =r$, $m + n \leq g$ and
\[ | \a_i, \a_j |= 1 = |\a_i, \b_i|, \textit{  and   } |\b_i, \b_j | =-1 \]
for all $i, j$.
\end{prop}

Such a group is said to be of \textbf{rank $r$} and \textbf{type $(m, n)$}. We will describe subgroups of $\J(\X)$ generated by half-periods using this property.  

Throughout this section we will make use of the list of groups for genus 3 curves as described in \cite{kyoto}.  

\subsection{Genus 3 curves with elliptic involutions}
Let $\X$ be a genus 3   curve and $\s \in \Aut (\X)$ an elliptic involution. Denote by $\pi$ the quotient map   $\pi : \X \to \X/ \< \s \>$.  We denote $\X / \< \s \>$ by $\E$.  Without loss of generality we assume that $\pi : \X \to \E$ is maximal. Then, $ \E \embd \J (\X)$.   The map $\pi :  \X  \to \E$ has four branch points. By picking the right origin
for $\E$, we can assume that the set of the branch points is the set of 2-torsion points on $\E$ which we denote by $\E [2]$.

The  points in the set   $\pi^\ast ( \E [2] ) \subset \J (\X)$, 
are called the \textbf{corresponding 2-torsion points of $\s$}. Accola and others have called them \textbf{derived half-periods of} $\s$. To simplify the notation we use 
\[\rmJ_\s := \pi^\ast (\E [2]).\]
Next, we give a more topological description of these points and the action of $\s$ on $H_1 (\X, \Z)$.  We
pick a homology basis $\{ A, B\}$ for $\E$ such that $A_2, A_3$ and $B_2, B_3$ are the lifting of
respectively $A$ and $ B$ in $H_1 (\X, \Z)$. Notice that $\s$ acts on $H_1 (\X, \Z)$ by
\[ \s (A_1, A_2, A_3, B_1, B_2, B_3) = (-A_1, A_3, A_2, -B_1, B_3, B_2)\]
Let $V^\s$ denote the $\s$-invariant subspace of $H_1 (\X,\Z)$ and
\[ \Phi : H_1 (\X, \Z) \to H_1 (\X, \Z / 2\Z)\]
the natural projection. Then we have the following:

\begin{lem}
The set of corresponding 2-torsion points  of $\s$ is the set $\Phi (V^\s)$.
\end{lem}

\proof The map $\pi^\ast : \E \to \J (\X)$ is given by
\[ \pi^\ast \begin{pmatrix} a \\ b \end{pmatrix} = \begin{pmatrix} 0 & a & a \\ 0 & b & b \end{pmatrix}, \]
see \cite[pg. 44]{Accola}. The point $\begin{pmatrix} a \\ b \end{pmatrix}$ is a 2-torsion point in $\E$. Hence,
$a, b \in (\frac 1 2 \Z/2\Z)$. Thus,
\[ \E [2]= \{\begin{pmatrix} 0 \\ 0 \end{pmatrix}, \ \begin{pmatrix} \frac 1 2 \\ 0 \end{pmatrix},  \
\begin{pmatrix} 0 \\ \frac 1 2  \end{pmatrix}, \begin{pmatrix} \frac 1 2 \\ \frac 1 2 \end{pmatrix} \]
Hence, the derived half-periods are
\[ \begin{pmatrix} 0 & \frac 1 2 & \frac 1 2 \\ 0 & 0 & 0 \end{pmatrix},
\, \, \begin{pmatrix} 0 & 0 & 0 \\ 0 & \frac 1 2  & \frac 1 2 \end{pmatrix},
\, \, \begin{pmatrix} 0 & \frac 1 2 & \frac 1 2 \\ 0 & \frac 1 2  & \frac 1 2 \end{pmatrix}, \]
which correspond to the $\s$-invariant subspace
\[ \{ B_2 +B_3, A_2 + A_3, A_2 + B_2 + A_3 + B_3 \} \]
This completes the proof.
\endproof
The proof of the following is intended in \cite{b-e-sh}. 
\begin{lem} Two  elliptic involutions   $\s, \t \in \Aut (\X)$ commute if and only if  $|\rmJ_\s \cap \rmJ_\t | =2$.
\end{lem}

%
%



We want to describe the properties of $\J (\X)$ in terms of 4-torsion elements of $\J (\X)$. Hence we also
define the following
\[ \JJ_\s : = \pi^\ast ( \E [4] ) \subset \J (\X),\]
Then $\JJ_\s $ is a group of order 16 in $\J (\X)$ and we have
\[ J_\s < \JJ_\s < \J (\X)\]
Hence, there are 12 quarter-periods (not including half-periods) which belong to $\s$.

\begin{prop}
Each point $\p \in \JJ_\s$ is a theta-null.
\end{prop}

\proof Let $\X$ be a genus 3 curve, $\O$ a period matrix of $\X$, $\Theta_\X$ its theta divisor, $\s \in \Aut
(\X)$ an elliptic involution, and $\pi : \X \to \X^\s$ the quotient map. Then, there is a half-period theta
null $\p \in J_s  $ such that for any $\a \in \X^\s$ we have \[\pi^\ast (\a) + \p \in \Theta_\X ,\] see
Accola, \cite[p. 88]{Accola}.  Hence, for any $\a \in \X^\s$ we have
\[\T (\p + \pi^\ast (\a), \O ) = 0\]
In particular, there are exactly 12  such points $\a \in \E[4] \setminus \E[2]$  such that  $\T (\p +
\pi^\ast (\a), \O ) = 0$. Since $\p$ has order 2 and all $\a_1, \dots , \a_{12}$ have order 4 then  all
points $\p + \pi^\ast (\a_i)$ have order 4,  for all $i=1, \dots , 12$. We denote all these quarter periods
$\f_i :=\p + \pi^\ast (\a_i)$. Thus, $\JJ_\s$ is the union of $ \{ \f_i \}_{i=1}^{12}$ and $J_\s$. Since,
points in $J_\s$ are theta-nulls and from the above all $\f_i$ are theta-nulls the conclusion holds.

\endproof

\noindent Notice that, $\E [4] \iso C_4 \times C_4$, it can be given as
\[ \E [4]=\left\{ \frac a b \, \, \left| \, \, a, b \in \left( \frac 1 4 \Z/ 4\Z \right)  \right.  \right\}  .\]
%
%
%
Next we intend to find necessary and sufficient conditions on half-periods and quarter-periods which will determine the automorphism group of a genus 3 curve.
  We summarize all the cases of non-hyperelliptic genus 3 curves in the following theorem. Notice that when we
say "there exist quarter periods theta-nulls" we always mean  "distinct" periods.

\begin{thm}\label{thm1} Let $\X$ be a genus 3 algebraic curve, $\O$ its period matrix, and $G=\Aut(\X)$ its group of
automorphisms. Then, the following hold:
\begin{enumerate}
\item     $C_2 \embd G$ if and only if there exist two    quarter periods theta-nulls $\f_1, \f_2 \in \J (\X)$ such that $\f_1 \neq \pm \f_2$,  $2 \f_1=2 \f_2 $ and  $|2 \f_1, \f_1+\f_2 |=1$.

\item if $V_4 \embd G$ then we have the following two cases:
\begin{enumerate}
\item  If $\X$ is hyperelliptic:  $V_4 \embd G$  if and only if then there are three   quarter periods
theta-nulls $\f_1, \f_2, \mathfrak c$ in $\J (\X)$ such that $2\f_1 = 2\f_2 = 2\mathfrak c$.

\item  If $\X$ is not hyperelliptic:  $V_4 \embd G$  if and only if     $\th (z, \O)$ vanishes to order 2 at one
half-period and there are two quarter periods $\f_1 , \f_2 \in \J (\X)$ such that  $\f_1 \neq \pm \f_2$,  $2\f_1 = 2\f_2 \neq 0$, and   $|2 \f_1, \f_1+\f_2 |=1$.
\end{enumerate}

\item    $C_3 \embd G$     if and only if there exist two    $\frac 1 6$--periods theta-nulls $\f_1, \f_2 \in \J (\X)$ such that
   $3 \f_1=3 \f_2 $ and  $|2 \f_1, \f_1+\f_2 |=1$.

\item If $G = C_2^3$ then $\th (z, \O)$ vanishes to order 2 at one half-period and vanishes to order one at
three quarter periods $\f_1, \f_2, \f_3 \in \J (\X)$ such that $2\f_1 = 2\f_2 = 2\f_3$.

\item     $S_3 \embd G$   if and only if there are four   quarter-periods theta-nulls $\f_1, \dots ,
\f_4 \in \J (\X)$ such that
\begin{itemize}
\item[i)] $2\f_1=2\f_2 \neq 2\f_3=2\f_4$

\item[ii)] $|2\f_1, \f_1+\f_2| = 1= |2\f_3, \f_3+\f_4|$

\item[iii)] $\< 2\f_1, \f_1+\f_2, 2\f_3, \f_3+\f_4 \>$ has type $(0, 2)$.
\end{itemize}

\item   $D_4 \embd G$   if and only if   there are four   quarter-periods theta-nulls $\f_1, \f_2, \f_3, \f_4 \in \J (\X)$ such that
\begin{itemize}
\item[i)] $2\f_1=2\f_2 \neq 2\f_3=2\f_4$

\item[ii)] $|2\f_1, \f_1+\f_2| = 1= |2\f_3, \f_3+\f_4|$

\item[iii)]  $ \< 2\f_1, \f_1+\f_2, 2\f_3, \f_3+\f_4 \>$  has type $(2, 1)$.
\end{itemize}

\item   $S_4 \embd G$ if and only if there are five   quarter-periods theta-nulls $\f_1, \f_2, \h_1, \h_2, \h_3 \in \Jac (\X)$ such that
\begin{itemize}
\item[i)] $2\f_1=2\f_2, 2\h_1=2\h_2=2\h_3$,

\item[ii)] $|2\f_1, \f_1+\f_2|=|2\h_1, \h_1+\h_2|=|2\h_1, \h_1+\h_3|=1$

\item[iii)] $\< 2 \f_1, \f_1 + \f_2, 2 \h_1, \h_1 + \h_2\> $ has type $(2,1)$

\item[iv)]  $\< 2 \f_1, \f_1 + \f_2, 2 \h_1, \h_1 + \h_3 \> $ has type $(0, 2)$
\end{itemize}

\item $C_4^2 \xs S_3   \embd G$ if and only if there are six   quarter-periods theta-nulls $\f_1, \dots , \f_6 \in \J (\X)$ such that
\begin{itemize}
\item[i)]   $2\f_1=2\f_2, 2\f_3=2\f_4, 2\f_5=2\f_6$

\item[ii)]    $|2\f_1, \f_1+\f_2|=|2\f_3, \f_3+\f_4|=|2\f_5, \f_5+\f_6 |=1$

\item[iii)] $H:=\< 2\f_1, \f_1+\f_2, 2\f_3, \f_3+\f_4\>$ has type $(2, 1)$.

\item[iv)] $\< 2\f_1, \f_1+\f_2, 2\f_3, \f_3+\f_4, 2\f_5, \f_5 + \f_6\>$ has rank 6.

\item[v)] if $y, z \in H$ satisfy $|x,y|=|x,z|=1$ for all $x \in H$, then $M:=\< y, z, 2\f_5, \f_5+\f_6\>$ has type $(2, 1)$.
\end{itemize}

\item     $  L_3(2) \embd G$ if and only if there are six   quarter-periods theta-nulls $\f_1, \dots , \f_6 \in \J (\X)$ such that they satisfy conditions $i) \dots iv) $ of Case 4) and $M$ has type $(0, 2)$.

%
%
%
%


%
%
%
%




\end{enumerate}
\end{thm}

\noindent The proof of the   theorem will take the rest of this paper.  Detailed proofs of some of the results here and generalizations to higher genus are intended in \cite{b-e-sh}.


Let us assume now that $\X$ is a genus 3 non-hyperelliptic curve such that $C_2 \embd \Aut (\X)$. Hence, $\X$ has an elliptic involution which we denote by $\s$. Thus there exists subgroups $\rmJ_\s, \JJ_\s$ in $\J (\X)$
such that these elements are theta-nulls.   We take
\[ \rmJ_\s= 
\left\{
\begin{pmatrix} 0 \\ 0 \end{pmatrix}, \ \begin{pmatrix} \frac 1 2 \\ 0 \end{pmatrix},  \
\begin{pmatrix} 0 \\ \frac 1 2  \end{pmatrix}, \begin{pmatrix} \frac 1 2 \\ \frac 1 2 \end{pmatrix} 
\right\}
\]
The elements $\E [4]$ are mapped to $H_1 (\X, \Z/ 4\Z)$ as follows
\[ \pi^\ast   \begin{pmatrix} a \\ b \end{pmatrix} =
\begin{pmatrix} \frac 1 2  & \frac a 2 & \frac a 2 \\ \frac 1 2 & \frac b 2 & \frac b 2 \end{pmatrix}.     \]
We take the elements
\[
\begin{array}{cc}
\f_1:=\begin{pmatrix} \frac 1 2  & \frac 1 4 & \frac 1 4 \\ \frac 1 2 & 0 & 0 \end{pmatrix}, &
\f_2:=\begin{pmatrix} \frac 1 2  & 0 &  0 \\ \frac 1 2 & \frac 1 4 & \frac 1 4 \end{pmatrix},
\end{array}
\]
Notice that $\f_1, \f_2$ are elements of order 4 and $\< \f_1\> \cap \<\f_2\> = \{ 0\}$. Hence,  $\JJ_\s \iso
\< \f_1\> \times \< \f_2\>$. Also
\[
\begin{array}{cc}
2\f_1:=\begin{pmatrix} 0  & \frac 1 2 & \frac 1 2 \\ \frac 1 2 & 0 & 0 \end{pmatrix}, &
2\f_2:=\begin{pmatrix} 0  & 0 &  0 \\ 0 & \frac 1 2 & \frac 1 2 \end{pmatrix},
\end{array}
\]
Hence,  $\rmJ_\s= \< 2\f_1, 2\f_2\>$ which is isomorphic to the Klein 4-group. Then, we have the following, part of which is proved in \cite[Cor. 5, pg. 53]{Accola}.

\begin{prop}\label{c2}
Let $\X$ be a genus 3 curve. Then, $\X$ has an elliptic involution $\s$ if and only if there exist  quarter
periods $\f_1, \f_2 \in  \J (\X)$ such that:

\begin{itemize}
\item[i)]  $H:=\< \f_1, \f_2 \> $ is isomorphic to $C_4 \times C_4$.

\item[ii)]  all elements of $H$ are theta-nulls

\item[iii)]  the subgroup $\< 2\f_1, 2\f_2\> \leq H$ is a subgroup of half-periods and isomorphic to the Klein
4-group.
\end{itemize}
\end{prop}



As noted above, we  denote the subgroup $\JJ_\s =\< \f_1, \f_2 \>$  (resp.,  $\rmJ_\s =\< 2\f_1, 2\f_2\>$) and call it the \textit{subgroup of the corresponding quarter-periods} (resp., corresponding half-periods) of $\s$. In addition to the  half-periods from 
\[ \rmJ_\s = \{ 0, 2\f_1, 2\f_2, 2(\f_1+\f_2)\},\]
there are exactly 12 quarter-periods in $\JJ_\s$.
%

Assume that   $\X$ is hyperelliptic and   $V_4 \embd G$.  Then $V_4$ has two elliptic involutions $\s, \t$ and the hyperelliptic involution $w=\s\t$. Then, $\s\t=\t\s$ and $| J_s \cap J_\t | =2$. Thus,  we can take
\[ \rmJ_\s = \left\{ 0, 2\f_1, 2\f_2, 2(\f_1+\f_2) \right\}, \quad \rmJ_\t = \left\{ 0, 2\h_1, 2\h_2, 2(\h_1+\h_2) \right\} \]
Thus, we have the following:
\begin{lem}
Let $\X$ be a genus 3 hyperelliptic curve and $G= \Aut (\X)$. Then, $V_4 \embd G$  if and only if   there are  quarter periods  $\f_1 , \f_2, \h_1, \h_2 \in \J (\X)$ such that
\begin{itemize}
\item[i)] the groups  $H_\f := \< \f_1, \f_2\> $ and $H_\h :=\< \h_1, \h_2\>$ are both isomorphic to $C_4 \times C_4$
\item[ii)] all elements of $H_\f$ and $H_\h$ are theta-nulls
\item[iii)] $H_\f \cap H_\h \iso C_2$.
\end{itemize}
\end{lem}
Now we assume that $\X$ is a non-hyperelliptic genus 3 curve and $V_4 \embd \Aut (\X)$.  Then, we have three elliptic involutions $\s, \t \in \Aut (\X)$ which all commute.   Thus, we have the following:
\begin{lem}Let $\X$ be a genus 3 non-hyperelliptic curve and $G= \Aut (\X)$.   $V_4 \embd G$  if and only if   there are  quarter-periods  $\f_1 , \f_2, \h_1, \h_2, \g_1, \g_2 \in \J (\X)$ such that
\begin{itemize}
\item[i)] $2\f_1 = 2\h_1=2\g_1$.

\item[ii)] the groups  $H_\f := \< \f_1, \f_2\> $, $H_\h :=\< \h_1, \h_2\>$, $H_\g :=\< \g_1, \g_2\>$ are all
isomorphic to $C_4 \times C_4$

\item[iii)] all elements of $H_\f, H_\h, H_\g$  are theta-nulls
\end{itemize}
\end{lem}

\noindent This case was also studied in \cite[Theorem 6, pg. 92]{Accola}.\\

 The automorphism group is $D_4$ and  $\X$ is non-hyperelliptic. We follow a more topological approach since there are five involutions in $D_4$ and it seems complicated to analyze the intersections among all the corresponding quarter-periods. Take $D_4$ as
\[D_4= \< \a, \b \ | \ \a^2=\b^4=1, \a\b\a=\b^3 \>\]
All involutions of $D_4$ are elliptic involutions. Also,  $Z (D_4) = \< \b^2\>$. Let $E:= \X / \< \b^2 \>$. Then, $\Aut (E) = D_4 / Z (D_4)$. Let
$\pi:   \Aut (\X)  \to \Aut (E)$ such that    $\a \to {\bar \a}$, $ \b  \to {\bar \b}$, and  $\a \b \to \overline{\a\b}$. 
Next, we find how $\bar \a, \bar \b, \overline{\a\b}$ act on $E$ and then lift them back to $\X$ to compute their action on homology. It is a simple exercise in covering spaces to determine that the action 
$ D_4 \times H_1 ( \X, \Z)  \to H_1 (\X, \Z)$  is given by
\begin{small}
\[
\begin{split}
\a: (A_1, \dots , B_3) & \to \left(A_1 - A_2 +A_3, -A_2, -A_3, B_1, -B_1-B_2, B_1-B_3 \right)\\
\b: (A_1, \dots , B_3) & \to (A_1 + A_2 -A_3, -A_1+A_2, A_1+A_3, B_1-B_2+B_3, \\
& B_1+B_3, -B_1+B_2)\\
\a \b : (A_1, \dots , B_3) & \to \left(A_1, A_1-A_2, -A_1-A_3, B_1+B_2+B_3, -B_2, -B_3 \right)\\
\end{split}
\]
\end{small}
Then, the action of all other involutions of $D_4$ is given by
\begin{small}
\[
\begin{split}
\a \b : (A_1, \dots , B_3) & \to (A_1, A_1-A_2, -A_1-A_3, B_1+B_2+B_3, -B_2, -B_3)\\
\b^2 : (A_1, \dots , B_3) & \to   (-A_1, A_3, A_2, -B_1, B_3, B_2)                \\
\a\b^2 : (A_1, \dots , B_3) & \to  (-A_1 + A_2 -A_3, -A_3, -A_2, -B_1, B_1-B_3, -B_1-B_2) \\
\a\b^3 : (A_1, \dots , B_3) & \to   (-A_1, -A_1-A_3, A_1-A_2 , -B_1-B_2+B_3, -B_3, -B_2) \\
\end{split}
\]
\end{small}
The invariant subspaces and their images in $H_1 (\X, \Z / 2\Z)$ are given below.
\begin{small}
\[
\begin{aligned}
& V^\a   =      \   \<B_1, 2A_1-A_2+A_3  \>,   
& \Phi(V^\a)= \left\{
\begin{bmatrix} \frac 1 2 , 0 , 0\\ 0, 0, 0 \end{bmatrix},
\begin{bmatrix} 0 , 0 , 0\\ 0, \frac 1 2, \frac 1 2 \end{bmatrix},
\begin{bmatrix} \frac 1 2 , 0 , 0\\ 0, \frac 1 2, \frac 1 2 \end{bmatrix}
\right\} & \\
& V^{\a \b}  =    \ \< A_1, 2B_1 +B_2-B_3\>,      
& \Phi(V^{\a \b})= \left\{
\begin{bmatrix} 0 , \frac 1 2 , \frac 1 2\\ 0, 0, 0 \end{bmatrix},
\begin{bmatrix} 0 , 0 , 0\\ \frac 1 2, 0, 0 \end{bmatrix},
\begin{bmatrix} 0 , \frac 1 2 , \frac 1 2\\ \frac 1 2, 0, 0 \end{bmatrix}
\right\}  & \\
& V^{\b^2}  =    \ \< A_2 - A_3, B_1+B_2-B_3\>,      
& \Phi(V^{\b^2})= \left\{
\begin{bmatrix} 0 , \frac 1 2 , \frac 1 2\\ 0, 0, 0 \end{bmatrix},
\begin{bmatrix} 0 , 0 , 0\\ 0, \frac 1 2, \frac 1 2 \end{bmatrix},
\begin{bmatrix} 0 , \frac 1 2 , \frac 1 2\\ 0, \frac 1 2, \frac 1 2 \end{bmatrix}
\right\}     &   \\
& V^{\a\b^2 } = \   \< A_2-A_3, B_1+B_2-B_3 \>,      
& \Phi(V^{\a\b^2})= \left\{
\begin{bmatrix} \frac 1 2 , \frac 1 2 , \frac 1 2\\ 0, 0, 0 \end{bmatrix},
\begin{bmatrix} 0 , 0 , 0\\ 0, \frac 1 2, \frac 1 2 \end{bmatrix},
\begin{bmatrix} \frac 1 2 , \frac 1 2 , \frac 1 2\\ 0, \frac 1 2, \frac 1 2 \end{bmatrix}
\right\}        \\
& V^{\a\b^3}  =\ \< A_1-A_2 + A_3, B_2-B_3 \>,       
& \Phi(V^{\a\b^3})= \left\{
\begin{bmatrix} 0 , \frac 1 2 , \frac 1 2\\ 0, 0, 0 \end{bmatrix},
\begin{bmatrix} 0 , 0 , 0\\ \frac 1 2, \frac 1 2, \frac 1 2 \end{bmatrix},
\begin{bmatrix} 0 , \frac 1 2 , \frac 1 2\\ \frac 1 2, \frac 1 2, \frac 1 2 \end{bmatrix}
\right\}   &    \\
\end{aligned}
\]
\end{small}
Then we have that $D_4 \embd G$   if and only if   there are four   quarter-periods theta-nulls $\f_1, \f_2, \f_3, \f_4 \in \J (\X)$ such that
\begin{itemize}
\item[i)] $2\f_1=2\f_2 \neq 2\f_3=2\f_4$
\item[ii)] $|2\f_1, \f_1+\f_2| = 1= |2\f_3, \f_3+\f_4|$
\item[iv)]  $ \< 2\f_1, \f_1+\f_2, 2\f_3, \f_3+\f_4 \>$  has type $(2, 1)$.
\end{itemize}


The automorphism group is  $S_4$. This is a subcase of case 2). Hence, there are four distinct quarter-periods theta-nulls $\f_1, \f_2, \h_1, \h_2 \in \J (\X)$  as in case 2). Notice that $D_4 \embd \Aut (\X)$ and there is another elliptic involution $\s \in \Aut (\X)$ such that $\s \not\in D_4$. Since $\s$ commutes with an involution from $D_4$ then there is a common derived half-period of $\s$ with the half-periods from above. So there exists a $\h_3$ such that $2\h_1=2\h_2=2\h_3$ and $\th (\h_3 )=0$.

Conversely, if there are five distinct quarter-periods theta-nulls $\f_1, \f_2, \h_1, \h_2, \h_3 \in \J (\X)$ as in Case 3)  then $D_4 \embd G$. The existence of a fifth quarter-period means that there is another
involution in $G$ which commutes with one of the involutions in $D_4$. Hence, $G$ is isomorphic to $S_4$ or the group with identity $(16, 13)$. The group $\< 2 \f_1, \f_1 + \f_2, 2 \h_1, \h_1 + \h_3 \> $ has type $(0, 2)$ which implies that $G \embd S_4$.



 
The automorphism group is  $C_4^2 \sem S_3$.  Then there is a dihedral group $D_4$ such that $D_4 < S_4  < C_4^2 \sem S_3 < G$. Hence, there exist $\f_1, \dots , \f_5 \in \J (\X)$ such
that they satisfy case 3).

There are involutions in $G$ which do not commute with some involution of $D_4$. Then, there exist another
half period $\f_6 \in \J (\X)$ which satisfies i) ... v).


The automorphism group is $G=GL_3(2)$.    Then,  $S_3$ and $D_4$ are both embedded in $G$. This locus is a sublocus of both $\M_3 (S_3)$ and $\M_3 (D_4)$. Hence, the involutions which come from the $S_3$ subgroup should generate a group of $(0, 2)$ type. The converse is similar.

The automorphism group is $S_3 \embd G$.  The proof of this case is similar to that of case 2). It is also proved in \cite{Wo} so we omit the details.

\subsection{Genus 3  curves with cyclic automorphism group, superelliptic curves.}

Let $\X$ be a non-hyperelliptic genus 3 curve with $\Aut (\X) = C_3$. Then there is a degree 3 covering $\pi
\colon \X \to \P_x^1$ branched at 5 points and with ramification index $3$ at each point. We take   the
branch points to be $\{ 0, 1, \infty, s, t\}$. We pick the points $P_{(0,0)}$ and $P_{(\infty, \infty)}$ to
be in the fibers of $0$ and $\infty$ respectively. Then the cover is $\pi (y) = y^3$ and the curve has
equation
\[ y^3 = x (x-1)(x-s)(x-t). \]
Let  $\e$ such that $\e^2+\e+1=0$. There is an order 3 automorphism $\s \in \Aut (\X)$ such that  $\s (x, y) \to (x, \e y)$. 
Let $\{ A_1, A_2, A_3, B_1, B_2, B_3\}$ be the canonical homology basis. Then
\[ \s (A_1, A_2, A_3, B_1, B_2, B_3) = (A_2, A_3, A_1, B_2, B_3, B_1)\]
We find the $\s$-invariant space of this action and then compute the thetanulls; details are intended in \cite{b-e-sh}.

\begin{prop}\label{C3}
Let $\X$ be a genus 3 algebraic curve, $\O$ its period matrix, and $G=\Aut(\X)$ its group of
automorphisms. Then,
 $C_3 \embd G$     if and only if there exist two    $\frac 1 6$--periods theta-nulls $\f_1, \f_2 \in \J (\X)$ such that
   $3 \f_1=3 \f_2 $ and  $|2 \f_1, \f_1+\f_2 |=1$.
\end{prop}


There are a few other groups which occur when $g=3$, but the corresponding loci for such groups has dimension zero and we do not discuss them here. 
The above cases finish the proof of the Theorem~\ref{thm1}. 

A similar project could be attempted for genus $g=4$ curves, since the list of groups and the corresponding among the loci is known; see \cite{nato-2} in this volume. 

\nocite{*}

\bibliographystyle{amsplain} 

\bibliography{theta}{}


\end{document}